\documentclass{gtpart}   
\usepackage{pinlabel}

\title[Erratum: Floer homology and splicing knot complements]
{Correction to the article:\\
Floer homology and splicing knot complements}%
\author{Eaman Eftekhary}%
\address{School of Mathematics, Institute for Research in Fundamental Sciences (IPM),
P. O. Box 19395-5746, Tehran, Iran}%
\email{eaman@ipm.ir}
\urladdr{http://math.ipm.ir/~eftekhary}



\arxivreference{arXiv:0802.2874}  
\arxivpassword{h2rkv}   

\usepackage{amsmath,amsthm,amsfonts,amscd,flafter,epsf,
amssymb,graphicx,color}
\usepackage{epsfig,manfnt}
\usepackage{amsfonts}
\usepackage{amssymb}
\usepackage{amsmath}
\usepackage{amsthm}
\usepackage{latexsym}
\usepackage{amscd}
\usepackage{flafter}
\usepackage{epsf}
\usepackage{tikz,tikz-cd}

\input{diagram}

\newcommand{\colvec}[2][.8]{%
  \scalebox{#1}{%
    \renewcommand{\arraystretch}{.8}%
    $\begin{pmatrix}#2\end{pmatrix}$%
  }
}
\newcommand{\lift}[2]{%
\setlength{\unitlength}{1pt}
\begin{picture}(0,0)(0,0)
\put(0,{#1}){\makebox(0,0)[b]{${#2}$}}
\end{picture}
}
\newcommand{\lowerarrow}[1]{%
\setlength{\unitlength}{0.03\DiagramCellWidth}
\begin{picture}(0,0)(0,0)
\qbezier(-28,-4)(0,-18)(28,-4)
\put(0,-14){\makebox(0,0)[t]{$\scriptstyle {#1}$}}
\put(28.6,-3.7){\vector(2,1){0}}
\end{picture}
}
\newcommand{\upperarrow}[1]{%
\setlength{\unitlength}{0.03\DiagramCellWidth}
\begin{picture}(0,0)(0,0)
\qbezier(-28,12)(0,25)(28,12)
\put(0,20){\makebox(0,0)[b]{$\scriptstyle {#1}$}}
\put(28.6,11.7){\vector(2,-1){0}}
\end{picture}
}
\newcommand{\strarrow}[1]{%
\setlength{\unitlength}{0.03\DiagramCellWidth}
\begin{picture}(0,0)(0,0)
\qbezier(-22,5)(0,5)(22,5)
\put(0,6){\makebox(0,0)[b]{$\scriptstyle {#1}$}}
\put(-22,4.8){\vector(-2,0){0}}
\end{picture}
}

\newcommand{\fbar}{\overline{f}}

\newcommand{\pphi}{\mathfrak{f}}
\newcommand{\pphibar}{\overline{\pphi}}
\newcommand{\CFDT}{\mathrm{CFD}}
\newcommand{\ovl}{\overline}
\newcommand{\ubar}{{\ovl{u}}}
\newcommand{\vbar}{{\ovl{v}}}
\newcommand{\wbarr}{{\ovl{w}}}

\newcommand{\Dd}{\mathfrak{D}}

\newcommand{\thetabar}{\overline{\theta}}

\newcommand{\B}{\mathbb{B}}

\newcommand{\HFT}{\mathrm{HF}}
\newcommand{\CFT}{\mathrm{CF}}
\newcommand{\HFKT}{\mathrm{HFK}}

\newcommand{\rank}{\mathrm{rnk}}

\newcommand{\ov}{\widehat}

\newcommand{\E}{\mathbb{E}}
\newcommand{\Sig}{\Sigma}
\newcommand{\ra}{\rightarrow}

\newcommand{\A}{\mathbb{A}}

\newcommand\Dual{\mathcal D}
\newcommand\Duality\Dual

\newcommand\relspinc{s}

\newcommand\x{\mathbf x}

\newcommand\y{\mathbf y}

\newcommand\ModSphere{\ModFlow\left({\mathbb S}\longrightarrow
\Sym^{g-1}(\Sigma_{1})\times \Sym^2(\Sigma_{2})\right)}
\newcommand\ModSpheres\ModSphere

\newcommand\UnparModSp{\widehat \ModSp}
\newcommand\UnparModFlow\UnparModSp

\newcommand\ModMaps{\mathcal M}
\newcommand\ModSp\ModMaps

\newcommand\alphas{\mbox{\boldmath$\alpha$}}

\newcommand\betas{\mbox{\boldmath$\beta$}}

\newcommand\spincrel\relspinc

\newtheorem{thm}{Theorem}[section]
\newtheorem{prop}[thm]{Proposition}
\newtheorem{cor}[thm]{Corollary}

\newtheorem{defn}[thm]{Definition}

\newtheorem{remark}[thm]{Remark}

\def\endproof{\relax\ifmmode\expandafter\endproofmath\else
  \unskip\nobreak\hfil\penalty50\hskip.75em\hbox{}\nobreak\hfil\bull
  {\parfillskip=0pt \finalhyphendemerits=0 \bigbreak}\fi}
\def\endproofmath$${\eqno\bull$$\bigbreak}
\def\bull{\vbox{\hrule\hbox{\vrule\kern3pt\vbox{\kern6pt}\kern3pt\vrule}\hrule}}

\newcommand{\Ker}{\mathrm{Ker}}

\newcommand{\Coker}{\mathrm{Coker}}

\newcommand{\ModSWfour}{\mathcal{M}}
\newcommand{\ModFlow}{\ModSWfour}

\newcommand\abuts\Rightarrow
\newcommand\Sym{\mathrm{Sym}}

\newcommand{\Hbb}{{\mathbb{H}}}
\newcommand{\Hinfinf}{{\Hbb_{\infty,\infty}}}
\newcommand{\Hinfone} {{\Hbb_{\infty,1}}}
\newcommand{\Honeinf} {{\Hbb_{1,\infty}}}
\newcommand{\Honeone} {{\Hbb_{1,1}}}
\newcommand{\Honezero} {{\Hbb_{1,0}}}
\newcommand{\Hzeroone} {{\Hbb_{0,1}}}
\newcommand{\Hzerozero} {{\Hbb_{0,0}}}

\newcommand{\Fbb}{\mathbb{F}}

\newcommand{\lra}{\longrightarrow}

\begin{document}

\begin{abstract}
This note corrects the mistakes in the splicing formulas of the paper
"Floer homology and splicing knot complements" \cite{Ef}. 
The mistakes are 
the result of the incorrect assumption that  for a knot $K$ inside a homology sphere $Y$, the involution on 
$\ov\HFKT(K)$ which corresponds to moving the basepoints by one 
full twist around $K$ is trivial. The correction implicitly involves 
considering the contribution from this (possibly non-trivial) involution
in a number of places. 
\end{abstract}
\maketitle

We first describe the incorrect assumption about the triviality of 
the involution on $\ov\HFKT(K)$ which corresponds to moving the basepoints by one full twist around a knot $K$ inside a homology sphere $Y$. The incorrect assumption
that this map is trivial first appears in Subsection 3.3 of 
"Floer homology and splicing knot complements" \cite{Ef}. It is then used
in Subsection 5.3 to give a simplified matrix presentation of the 
splicing formula. For the trefoils, the involution map happens to 
be trivial, and thus the computations of Subsection 6.1 remain 
unchanged. Some modifications are necessary to the computations
of subsection 6.2. Finally, there are some changes that should 
be made to the matrix presentations in the appendix.
The corrected versions of subsections 3.3, 5.3 and 6.2, as well as 
the corrected form of the appendix appear in this note.\\
 
The original form of the splicing formulas, which are based on the aforementioned wrong assumption, are in fact not incorrect. 
However,  the arguments of this paper are not sufficient for 
showing this and the corrections of this note are thus very crucial.
The issue is further discussed in the sequel 
\cite{Ef-essential}.\\

The author would like to thank Adam Levine for pointing out the mistake, and for some very helpful discussions. 

\subsection*{The involution of $\ov\HFKT$}
In Subsection 3.3 of \cite{Ef}, the author considers the isomorphisms 
$\tau_\bullet:\Hbb_\bullet(K)\ra \Hbb_\bullet(K)$ for $\bullet\in\{0,1,\infty\}$ coming from 
changing the role of the two punctures $z_1$ and $z_2$ in the Heegaard diagram $(\Sig,\alphas,\betas;z_1,z_2)$. 
It was claimed that
$\tau_\bullet^2$ is the identity map. However, the isomorphism 
$\tau_\bullet^2$
is in fact the involution map $\varsigma_\bullet$
studied by Sarkar~\cite{Sarkar} and by Hendricks and Manolescu~\cite{HM}, which is not necessarily trivial.
This map corresponds to moving the punctures $z_1$ and $z_2$ one full twist along the knot $K$ and its square is the identity map.
The map $\varsigma_\infty$ is in general non-trivial, e.g. for 
the figure-eight knot.
Correspondingly, it is not {\emph{a priori}} true that 
\[\overline{\mathfrak{f}}_0=\tau_\infty\circ \mathfrak{f}_0\circ \tau_1,\ \ 
\overline{\mathfrak{f}}_1=\tau_0\circ \mathfrak{f}_1\circ \tau_\infty\ \ \text{and}\ \ 
\overline{\mathfrak{f}}_\infty=\tau_1\circ \mathfrak{f}_\infty\circ \tau_0.  
\]  
Without the above inequalities, many of the matrix computations in Subsection 5.3 need significant changes and the computation of 
Subsection 6.3 should also be modified. Some changes are also 
necessary in the later parts of the appendix.
The revised forms of these parts follow.\\

Nevertheless, the $3$ equalities
\[\overline{\mathfrak{f}}_0=\tau_\infty\circ \mathfrak{f}_0\circ \tau_1,\ \ 
\overline{\mathfrak{f}}_1=\tau_0\circ \mathfrak{f}_1\circ \tau_\infty\ \ \text{and}\ \ 
\overline{\mathfrak{f}}_\infty=\tau_1\circ \mathfrak{f}_\infty\circ \tau_0.  
\]  
are true for non-trivial reasons that will not be discussed here and are not relevant for the purposes of this paper. A detailed discussion of this  claim appears in \cite{Ef-essential}.

\setcounter{section}{3}  
\setcounter{subsection}{3}
\setcounter{equation}{3}
\subsection*{Subsection $3.3'$: Some properties of the maps
$\pphi_\bullet(K)$ and $\pphibar_\bullet(K)$}
Our first observation is that changing the orientation of the knot $K$, and
 correspondingly that of $K_{1}$ and $K_{0}$, corresponds to changing  
 the markings $u,v,w$ with $\ubar,\vbar,\wbarr$ in
 Figure 1.  Suppose that $(\Sig,\alphas,\betas;z_1,z_2)$ represents
$K_\bullet$, meaning that an oriented longitude for $K_\bullet$ is constructed from gluing
an oriented arc on $\Sig$ from $z_1$ to $z_2$ in the complement of $\alphas$ and an 
oriented arc on $\Sig$ from $z_2$ to $z_1$ in the complement of $\betas$. Then 
 $(\Sig,\alphas,\betas;z_2,z_1)$ is a Heegaard diagram for $-K_\bullet$
(the knot $K_\bullet$ with the reverse orientation) while  
 $(-\Sig,\betas,\alphas;z_2,z_1)$ is a Heegaard diagram for  $K_\bullet$.
  The chain complexes associated with the above three 
Heegaard diagrams are identical. 
Heegaard moves give chain homotopy equivalences
\begin{displaymath}
\tau_\bullet(K) \co \ov{\CFT}(\Sig,\alphas,\betas;z_1,z_2)=
\ov{\CFT}(-\Sig,\betas,\alphas;z_2,z_1)\lra 
\ov{\CFT}(\Sig,\alphas,\betas;z_1,z_2).
\end{displaymath}
These chain homotopy equivalences induce the isomorphisms
\begin{displaymath}
\tau_\bullet(K) \co \Hbb_\bullet(K)\lra \Hbb_\bullet(K),\ \ \ \bullet\in\{0,1,\infty\}.
\end{displaymath}
In terms of these isomorphisms  
\begin{equation}\label{eq:duality}
\begin{split}
&\pphibar_0(K)=\tau_{\infty}(K)\circ \pphi_0(K)\circ \tau_{1}(K)^{-1},\\ 
&\pphibar_1(K)=\tau_{0}(K)\circ \pphi_1(K)\circ \tau_{\infty}(K)^{-1}\ \  
 \text{and}\\
&\pphibar_\infty(K)=\tau_{1}(K)\circ \pphi_\infty(K)\circ 
\tau_{0}(K)^{-1}.
\end{split}
\end{equation}
Note however, that the equality
$\thetabar(K)=\tau_\infty(K)\theta(K)\tau_0(K)^{-1}$ 
is only satisfied for the induced maps 
from $\Ker(\pphibar_\infty(K))$ to $\Coker(\pphibar_0(K))$.
\\

The exactness of the sequences in (3) implies that  in appropriate 
decompositions 
\begin{equation}\label{eq:decompositions}
\begin{split}
&\Hbb_0(K)=\frac{\Hbb_0(K)}{\Ker(\pphi_\infty(K))}\oplus \Ker(\pphi_\infty(K))
=:\A_\infty(K)\oplus \A_1(K),\\ 
&\Hbb_1(K)=\frac{\Hbb_1(K)}{\Ker(\pphi_0(K))}\oplus \Ker(\pphi_0(K))
=:\A_0(K)\oplus \A_\infty(K)\ \ \ \text{and}\\
 &\Hbb_\infty(K)=\frac{\Hbb_\infty(K)}{\Ker(\pphi_1(K))}\oplus \Ker(\pphi_1(K))
=:\A_1(K)\oplus \A_0(K)
\end{split}
\end{equation}
we have $\pphi_\bullet(K)=\colvec{0&0\\ I_{a_\bullet(K)}&0}$,
where $a_\bullet(K)$ denotes the rank of $\A_\bullet(K)$ for every $\bullet\in \{0,1,\infty\}$.
In this basis we may present the matrices $\tau_\bullet(K)$  as
\begin{displaymath}
\tau_\bullet(K)=\left(\begin{array}{cc}
A_\bullet(K)& B_\bullet(K)\\ C_\bullet(K)& D_\bullet(K)
\end{array}\right),\ \ \ \
\bullet\in\{0,1,\infty\}.
\end{displaymath}

The map $B_0(K)$ corresponds to the induced map
$$\tau_0(K) \co \Ker(\pphi_\infty(K))\ra \frac{\Hbb_0(K)}{\Ker(\pphi_\infty(K))}.$$
The decomposition
$\Hbb_0(K)=\A_\infty(K)\oplus \A_1(K)$
may be modified using a change of basis  of the form
$P_X=\colvec{
I&0\\ -X&I}$, which does not change the block presentations of the 
maps $\pphi_\infty(K)$ and $\pphi_1(K)$. In the new 
basis $\tau_0(K)$ has the following presentation:
\begin{displaymath}
\begin{split}
\tau_0(K)&=\left(\begin{array}{cc}
I&0\\ -X&I
\end{array}\right)\left(\begin{array}{cc}
A_0(K)&B_0(K)\\
C_0(K)&D_0(K)
\end{array}\right)\left(\begin{array}{cc}
I&0\\ -X&I
\end{array}\right)
\\&
=\left(\begin{array}{cc}
A_0(K)-B_0(K)X&B_0(K)\\
\star & -XB_0(K)+D_0(K)
\end{array}\right)
\end{split}
\end{displaymath}
If $B_0(K)$ is injective, we may thus assume that $D_0(K)=0$, while if 
$B_0(K)$ is surjective, we may assume that $A_0(K)=0$.
With a similar reasoning, if $B_\bullet(K)$ is injective we may assume that $D_\bullet(K)=0$,
while if $B_\bullet(K)$ is surjective we may assume that $A_\bullet(K)=0$. \\

In the above decompositions for $\Hbb_\bullet(K)$, the map 
$\theta(K):\Hbb_0(K)\ra \Hbb_\infty(K)$ takes the form
\begin{displaymath}
\theta(K)=\left(\begin{array}{cc}
X&I\\Z&Y
\end{array}\right),
\end{displaymath}
since the induced map from $\A_1(K)\subset \Hbb_0(K)$ to 
$\A_1(K)\subset \Hbb_\infty(K)$ is the inverse of the map induced by 
$\pphi_1(K)$, i.e. the identity. Moreover, since the rank of $\theta(K)$ is the 
same as the rank of $\pphi_1(K)$, we conclude that $Z=YX$. Applying the 
change of basis $P_Y$ on $\Hbb_0(K)$ and the corresponding change of basis 
$P_X$ on $\Hbb_\infty(K)$, $\theta(K)$ takes the form
\begin{displaymath}
\left(\begin{array}{cc}I&0\\ -Y&I\end{array}\right)
\left(\begin{array}{cc}X&I\\ YX&Y\end{array}\right)
\left(\begin{array}{cc}I&0\\ -X&I\end{array}\right)=
\left(\begin{array}{cc}0&I\\ 0&0\end{array}\right).
\end{displaymath}
It is thus possible to choose the above decompositions so that 
$\theta(K)=\left(\begin{array}{cc}0&I\\ 0&0\end{array}\right)$.
If this is the case, the $2\times 2$ presentation of 
$\tau_\infty(K)^{-1}\thetabar(K)\tau_0(K)$ would be of the form
\begin{displaymath}
\tau_\infty(K)\thetabar(K)\tau_0(K)=\left(\begin{array}{cc}
M&I\\Q&P
\end{array}\right),
\end{displaymath}
and since the ranks of $\theta(K)$ and $\thetabar(K)$ are the same, we find
$Q=PM$.

\setcounter{thm}{4}
\setcounter{section}{5}  
\setcounter{subsection}{3}
\setcounter{equation}{17} 
\subsection*{Subsection $5.3'$. Simplifications of the splicing formula}\label{subsec:simplifications}
We now apply Lemma 2.3 to the splicing formula of  
Proposition 5.3 and make some cancellations. The first cancellation 
comes from setting $C=\Hbb$, $A=\Honeone$ and 
$$B=\left(\Hinfinf\oplus \Honeinf\oplus \Hinfone\right)\oplus
\left( \Hzeroone\oplus \Honezero\oplus \Hzerozero\right)=\E_1\oplus \E_2.$$
We thus have 
$\ov\HFT (Y)=H_*(B,d_B)$, where
\begin{displaymath}
d_B=\colvec{
0&\pphi_0^1\otimes I&I\otimes \pphi_0^2&\theta^1\otimes \pphibar_0^2&
\pphibar_0^1\otimes \theta^2&\Gamma\\
&&&\\
0&0&0&\Phi&I\otimes (\pphi_0^2\circ \pphibar_\infty^2)&
\pphi_\infty^1\otimes\thetabar^2\\
0&0&0&(\pphi_0^1\circ \pphibar_\infty^1)\otimes I&\Psi &
\thetabar^1\otimes \pphi_\infty^2\\
0&0&0&0&0&I\otimes \pphibar_\infty^2\\
0&0&0&0&0&\pphibar_\infty^1\otimes I\\
0&0&0&0&0&0\\
},\end{displaymath}
with 
$\Gamma=\pphibar_0^1\circ \pphi_\infty^1)\otimes
\pphibar_0^2\circ \pphi_\infty^2+
\theta^1\otimes \thetabar^2+\thetabar^1\otimes \theta^1$,
$\Phi=\pphibar_\infty^1\otimes \pphi_0^2+\pphi_\infty^1\otimes \pphibar_0^2$ and 
$\Psi=\pphi_0^1\otimes \pphibar_\infty^2+\pphibar_0^1\otimes \pphi_\infty^2$.\\

The dimension of the $\Fbb$-vector space $H_*(B,d_B)$ only depends on 
the rank of the kernel and the cokernel of the matrix $d_B$. Define a pair of matrices 
$M_1$ and $M_2$ equivalent if $\Ker(M_1)\simeq \Ker(M_2)$ and 
$\Coker(M_1)\simeq \Coker(M_2)$. For a matrix $M$ let 
$$\imath(M):=\Ker(M)\oplus \Coker(M)\ \ \ \text{and}\ \ \ i(M):=\rank(\imath(M)).$$
If $M_1$ and $M_2$ are equivalent matrices then $\imath(M_1)\simeq \imath(M_2)$
and $i(M_1)=i(M_2)$.\\
 
We make a change of basis for $\E_{2}$ which is given by the matrix 
\begin{displaymath}
\colvec{\tau_0(K_1)\otimes \tau_1(K_2)&0&0\\
0&\tau_1(K_1)\otimes \tau_0(K_2) &0\\
0&0& \tau_0(K_1)\otimes \tau_0(K_2)
}.\end{displaymath}
The matrix $d_B$ is thus equivalent to the matrix
\begin{displaymath}
d_B'=\colvec{
0&\pphi_0^1\otimes I&I\otimes \pphi_0^2&
\theta^1\tau_0^1\otimes \tau^2_\infty\pphi_0^2
&\tau^1_\infty\pphi_0^1\otimes\theta^2\tau_0^2
&\Gamma
\\
0&0&0&\Phi
&\tau_1^1\otimes \pphi_0^2\tau_1^2\pphi_\infty^2 &
\pphi_\infty^1\tau_0^1\otimes\thetabar^2\tau_0^2\\
0&0&0&\pphi_0^1\tau_1^1\pphi_\infty^1\otimes \tau_1^1 &
\Psi
&\thetabar^1\tau_0^1\otimes \pphi_\infty^2\tau_0^2\\
0&0&0&0&0&I\otimes \pphi_\infty^2\\
0&0&0&0&0&\pphi_\infty^1\otimes I\\
0&0&0&0&0&0\\
}.
\end{displaymath}
\begin{align*}
\Phi&= \tau_1^1\pphi_\infty^1\otimes \pphi_0^2\tau_1^2+
\pphi_\infty^1\tau_0^1\otimes \tau_\infty^2\pphi_0^2,\\
\Psi&= \pphi_0^1\tau_1^1\otimes \tau_1^2\pphi_\infty^2+
\tau_\infty^1\pphi_0^1\otimes \pphi_\infty^2\tau_0^2\ \ \ \text{and}\\
\Gamma&= \pphibar_0^1\pphi_\infty^1\tau_0^1\otimes 
\pphibar_0^2\pphi_\infty^2\tau_0^2
+\theta^1\tau_0^1\otimes \thetabar^2\tau_0^2
+\thetabar^1\tau_0^1\otimes \theta^2\tau_0^2.
\end{align*}
Let us use the decompositions of (\ref{eq:decompositions}) for $K_1$ and $K_2$
to obtain  a $24\times 24$ block decomposition of $d_B'$. 
Moreover, following the discussion at the end of Subsection 3.3
we may assume that in the corresponding decompositions,
$$\theta^i=\left(\begin{array}{cc}
0&I\\ 0&0
\end{array}\right) \ \ \ \text{and}\ \ \ 
(\tau_\infty^i)^{-1}\thetabar^i\tau_0^i=\left(\begin{array}{cc}
M^i&I\\ P^iM^i&P^i
\end{array}\right).$$
Each entry in the above $6\times 6$ decomposition for $d_\B'$ 
corresponds to a $4\times 4$ sub-matrix of the aforementioned $24\times 24$
decomposition. For instance the $(1,4)$ entry 
$\theta^1\tau_0^1\otimes \tau_\infty^2\pphi_0^2$ 
corresponds to 
\begin{displaymath}
\colvec{0&I\\ 0&0}\colvec{A_0^1&B_0^1\\ C_0^1& D_0^1}\otimes 
\colvec{A_\infty^2& B_\infty^2\\ C_\infty^2& D_\infty^2}\colvec{0&0\\ I&0}=
\colvec{C_0^1& D_0^1\\0&0}\otimes \colvec{B_\infty^2&0\\ D_\infty^2&0}=
\colvec{
C_0^1\otimes B_\infty^2& 0& D_0^1\otimes B_\infty^2&0\\
C_0^1\otimes D_\infty^2& 0& D_0^1\otimes D_\infty^2 &0\\
0&0&0&0\\
0&0&0&0
}.
\end{displaymath}
For another instance, note that the $(3,5)$ entry  corresponds to 
\begin{align*}
\colvec{0&0&0&0\\
0&0&0&0\\
A_1^1\otimes B_1^2&0 & B_1^1\otimes B_1^2&0\\
A_1^1\otimes D_1^2&0 & B_1^1\otimes D_1^2&0 }&+
\colvec{0 &0 &0&0\\
B_\infty^1\otimes A_0^2&B_\infty^1\otimes B_0^2&0&0\\
0 &0 &0&0\\
D_\infty^1\otimes A_0^2&D_\infty^1\otimes B_0^2&0&0
}\\&=
\colvec{0&0&0&0\\
B_\infty^1\otimes A_0^2&B_\infty^1\otimes B_0^2&0&0\\
A_1^1\otimes B_1^2&0 & B_1^1\otimes B_1^2&0\\
&&&\\
A_1^1\otimes D_1^2+D_\infty^1\otimes A_0^2&
D_\infty^1\otimes B_0^2 & B_1^1\otimes D_1^2&0 
}.
\end{align*}
The aforementioned $24\times 24$ decomposition
 includes identity matrices as the entries determined by the following block coordinates:
\begin{align*}
(2,9), \ \ (3,5),\  \ (4,6), \ \  (14,21),\ \  (16,23)\ \ \text{and}\ \ (20,22).
\end{align*}
We use the above $6$ identity matrices for cancellation to obtain an equivalent 
matrix $d=\left(\begin{array}{cc}
0& D\\ 0&0\end{array}\right)$ 
over $\B_1\oplus \B_2$, where 
 $\A_{\bullet \star}=\A_\bullet(K_1)\otimes \A_\star(K_2)$ and
\begin{displaymath}
\begin{split}
&\B_1=\A_{11}\oplus \A_{\infty 1}\oplus \A_{\infty 0}\oplus \A_{1 \infty}\oplus
\A_{0 \infty}\oplus \A_{0 0}\ \ \text{and}\\
&\B_2=\A_{\infty 0}\oplus \A_{1 0}\oplus \A_{\infty \infty}\oplus \A_{0 \infty}\oplus
\A_{0 1}\oplus \A_{1 1}.
\end{split}
\end{displaymath}
Re-arrange the rows and the columns of the matrix $D$, so that $D$ corresponds to the 
rows $11,7,8,10,12,1$ and the columns $19,13,15,17,18,24$ in the above 
$24\times 24$ decomposition to obtain the  following matrix
\begin{displaymath}
\colvec{
B_1^1\otimes B_1^2&B_1^1\otimes A_1^2&0&A_1^1\otimes B_1^2&0&0\\
0&A_0^1\otimes B_\infty^2 &B_0^1\otimes B_\infty^2 &0&0&B_0^1\otimes 
(A_\infty^2+ B_\infty^2 P^2)\\
D_1^1\otimes B_1^2& D_1^1\otimes A_1^2
+A_0^1\otimes D_\infty^2
& B_0^1\otimes D_\infty^2& C_1^1\otimes B_1^2&0&
B_0^1\otimes (C_\infty^2+ D_\infty^2 P^2)\\
0&0&0& B_\infty^1\otimes A_0^2   &   B_\infty^1\otimes B_0^2  
 &(A_\infty^1+B_\infty^1 P^1)\otimes B_0^2 \\
B_1^1\otimes D_1^2&B_1^1\otimes C_1^2&0  &
D_\infty^1\otimes A_0^2
+A_1^1\otimes D_1^2 & D_\infty^1\otimes B_0^2
&(C_\infty^1+ D_\infty^1 P^1) \otimes B_0^2\\
0&C_0^1\otimes B_\infty^2& D_0^1\otimes B_\infty^2&
B_\infty^1\otimes C_0^2&B_\infty^1\otimes D_0^2& 
\Gamma}
\end{displaymath}
with
\begin{align*}
\Gamma=
B_\infty^1 \ovl{B}_1^1 B_0^1\otimes B_\infty^2 \ovl{B}_1^2 B_0^2+
(A_\infty^1+ B_\infty^1 P^1)\otimes D_0^2+
D_0^1\otimes (A_\infty^2+ B_\infty^2 P^2).
\end{align*}
Here, we assume that 
\begin{align*}
\left(\tau_\bullet^i\right)^{-1}=\left(\begin{array}{cc}
\ovl{A}_\bullet^i&\ovl{B}_\bullet^i\\
\ovl{C}_\bullet^i&\ovl{D}_\bullet^i
\end{array}\right),\ \ \ \ \bullet\in\{0,1,\infty\},\ i=1,2,
\end{align*}
are the corresponding block decompositions for the inverses of 
$\tau_\bullet^i$.

This matrix is in turn equivalent to the matrix $\Dd=\Dd(K_1,K_2)$ below,
which is obtained by adding $I\otimes P^2$ times the third column and $P^1\otimes I$
times the fifth column to the last column of the above matrix:
\begin{displaymath}
\Dd=\colvec{
B_1^1\otimes B_1^2&C_1^1\otimes A_1^2&0&A_1^1\otimes B_1^2&0&0\\
0&A_0^1\otimes B_\infty^2 &B_0^1\otimes B_\infty^2 &0&0&B_0^1\otimes A_\infty^2\\
&&&&&\\
D_1^1\otimes B_1^2&
D_1^1\otimes A_1^2
+A_0^1\otimes D_\infty^2
& B_0^1\otimes D_\infty^2& C_1^1\otimes B_1^2&0&B_0^1\otimes C_\infty^2\\
&&&&&\\
0&0&0& B_\infty^1\otimes A_0^2   &   B_\infty^1\otimes B_0^2  
 &A_\infty^1\otimes B_0^2 \\
&&&&&\\
B_1^1\otimes D_1^2&B_1^1\otimes C_1^2&0  &
D_\infty^1\otimes A_0^2
+A_1^1\otimes D_1^2
& D_\infty^1\otimes B_0^2
&C_\infty^1 \otimes B_0^2\\
&&&&&\\
0&C_0^1\otimes B_\infty^2& D_0^1\otimes B_\infty^2&
B_\infty^1\otimes C_0^2&B_\infty^1\otimes D_0^2& 
\Psi},
\end{displaymath}
where 
$\Psi=A_\infty^1\otimes D_0^2+D_0^1\otimes A_\infty^2
+X^1\otimes X^2$ and 
$X^i=X(K_i)=B_\infty^i \ovl{B}_1^i B_0^i$ for $i=1,2$.
\\

Combining Proposition 5.3 with the above observations we find:
\begin{prop}\label{prop:splicing-2}
Let $K_i\subset Y_i$ denote null-homologous knots for $i=1,2$ and $Y(K_1,K_2)$ 
denote the three-manifold obtained by splicing the complement of $K_1$
with the complement of $K_2$. With the above definition of $\Dd(K_1,K_2)$
\begin{displaymath}
\begin{split}
\ov{\HFT}\left(Y(K_1,K_2),\Fbb\right)
&\simeq \imath\left(\Dd(K_1,K_2)\right)\\
\end{split}
\end{displaymath}
\end{prop}

\begin{cor}\label{cor:theta-independence}
The splicing formula of Proposition 5.3 is independent 
of the choice of extensions $\theta^i$ and $\thetabar^i$.
\end{cor}
\begin{proof}
The fact that the matrices $P^i$ and $M^i$ do not appear in the matrix 
$\Dd(K_1,K_2)$ implies that the choice of the extensions 
$\theta^i,\thetabar^i \co \Hbb_0^i\ra \Hbb_\infty^i$ does not change the rank of the 
homology group in the splicing formula of 
Proposition~5.3 or
Theorem~1.1.
\end{proof}
With the above corollary in place, the proof of Theorem~1.1 is 
now complete.

\begin{defn}
 For a pair of knots $K_i\subset Y_i$, $i=1,2$ define 
\begin{displaymath}
\begin{split}
\chi(K_1,K_2):&=\Big(h_1(K_1)-h_\infty(K_1)\Big)\Big(h_1(K_2)-h_\infty(K_2)\Big)\\
&\ \ \ \ \ \ \ \ \ \ \ \ 
-\Big(h_1(K_1)-h_0(K_1)\Big)\Big(h_1(K_2)-h_0(K_2)\Big)\\
\end{split}
\end{displaymath}
\end{defn}
Note that $\chi(K_1,K_2)$ is in fact the difference between the ranks of $\B_1=\B_1(K_1,K_2)$
and $\B_2=\B_2(K_1,K_2)$. In the corresponding $\Z/2\Z$-grading on $\B_1\oplus \B_2$,
$\chi(K_1,K_2)$ is thus the Euler characteristic of the chain 
complex $(\B_1\oplus \B_2,d)$. 
\begin{cor}\label{cor:bound-on-rank}
With the above notation fixed, 
\begin{displaymath}
\rank\left(\ov\HFT(Y(K_1,K_2))\right)\geq 
\left|\chi(K_1,K_2)\right|.
\end{displaymath}
\end{cor}
\begin{proof}
It is enough to note that
\begin{displaymath}
\chi(K_1,K_2)
=\rank\big(\Ker(\Dd(K_1,K_2))\big)-\rank\big(\Coker(\Dd(K_1,K_2))\big).
\end{displaymath}
\end{proof}

Consider the matrices
\begin{displaymath}
P_L=\colvec{
I\otimes {\ovl{A}}_1^2&0&0&0&I\otimes {\ovl{B}}_1^2&0\\
0&I\otimes {\ovl{A}}_\infty^2&I\otimes {\ovl{B}}_\infty^2&0&0&0\\
0&I\otimes {\ovl{C}}_\infty^2&I\otimes {\ovl{D}}_\infty^2&0&0&0\\
0&0&0&I\otimes {\ovl{A}}_0^2&0&I\otimes {\ovl{B}}_0^2\\
I\otimes {\ovl{C}}_1^2&0&0&0&I\otimes {\ovl{D}}_1^2&0\\
0&0&0&I\otimes {\ovl{C}}_0^2&0&I\otimes {\ovl{D}}_0^2\\
}
\end{displaymath}
and 
\begin{displaymath}
P_R=\colvec{
{\ovl{D}}_1^1\otimes I&0&0&{\ovl{C}}_1^1\otimes I&0&0\\
0&{\ovl{A}}_0^1\otimes I&{\ovl{B}}_0^1\otimes I&0&0&0\\
0&{\ovl{C}}_0^1\otimes I&{\ovl{D}}_0^1\otimes I&0&0&0\\
{\ovl{B}}_1^1\otimes I&0&0&{\ovl{A}}_1^1\otimes I&0&0\\
0&0&0&0&{\ovl{D}}_\infty^1\otimes I& {\ovl{C}}_\infty^1\otimes I\\
0&0&0&0&{\ovl{B}}_\infty^1\otimes I& {\ovl{A}}_\infty^1\otimes I\\
}.
\end{displaymath}
Both $P_R$ and $P_L$ are invertible and
 $\Dd(K_1,K_2)$ is  equivalent to  $\Dd'(K_1,K_2)=P_L\Dd(K_1,K_2) P_R$.
The matrix $\Dd'(K_1,K_2)$ has the following block presentation.
\begin{displaymath}
\colvec[.7]{
D_\infty^1\ovl{B}_1^1\otimes \ovl{B}_1^2A_0^2&
B_1^1\ovl{A}_0^1\otimes I&
B_1^1\ovl{B}_0^1\otimes I&
D_\infty^1 {\ovl{A}}_1^1\otimes \ovl{B}_1^2A_0^2&
I\otimes \ovl{B}_1^2B_0^2&0\\
I\otimes \ovl{B}_\infty^2 B_1^2 &
D_1^1\ovl{A}_0^1\otimes \ovl{B}_\infty^2 A_1^2 &
D_1^1\ovl{B}_0^1\otimes \ovl{B}_\infty^2 A_1^2& 0&
B_0^1\ovl{B}_\infty^1\otimes I&
B_0^1 \ovl{A}_\infty^1\otimes I\\
I\otimes \ovl{D}_\infty^2 B_1^2& \Psi_1& D_1^1 \ovl{B}_0^1\otimes \ovl{D}_\infty^2 A_1^2 
&0&0&0\\
B_\infty^1 \ovl{B}_1^1\otimes I&0&I\otimes \ovl{B}_0^2B_\infty^2& 
B_\infty^1 \ovl{A}_1^1\otimes I
&\Gamma_1
&\Gamma_2\\
D_\infty^1 \ovl{B}_1^1\otimes \ovl{D}_1^2A_0^2&0&0&\Psi_2&I\otimes \ovl{D}_1^2B_0^2&0\\
0&0&I\otimes \ovl{D}_0^2B_\infty^2 &0&\Gamma_3&\Gamma_4
}
\end{displaymath}
with
\begin{align*}
\Psi_1&=I\otimes I+D_1^1 \ovl{A}_0^1\otimes \ovl{D}_\infty^2 A_1^2\\
\Psi_2&=I\otimes I+D_\infty^1 \ovl{A}_1^1\otimes \ovl{D}_1^2 A_0^2\\
\Gamma_1&=D_0^1 \ovl{B}_\infty^1\otimes \ovl{B}_0^2A_\infty^2
+X^1\ovl{B}_\infty^1\otimes \ovl{B}_0^2X^2\\
\Gamma_2&=D_0^1 \ovl{A}_\infty^1\otimes \ovl{B}_0^2A_\infty^2
+X^1\ovl{A}_\infty^1\otimes \ovl{B}_0^2X^2\\
\Gamma_3&=D_0^1 \ovl{B}_\infty^1\otimes \ovl{D}_0^2A_\infty^2
+X^1\ovl{B}_\infty^1\otimes \ovl{D}_0^2X^2\\
\Gamma_4&=I\otimes I+D_0^1\ovl{A}_\infty^1\otimes \ovl{D}_0^2 A_\infty^2
+X^1 \ovl{A}_\infty^1\otimes \ovl{D}_0^2X^2,
\end{align*}
and is easier to use in actual computations. Note that 
$$\imath(\Dd'(K_1,K_2))\simeq \imath (\Dd(K_1,K_2))\simeq \ov\HFT(Y(K_1,K_2),\Fbb).$$

\setcounter{thm}{0}
\setcounter{subsection}{2}
\setcounter{section}{6}  
\setcounter{equation}{21} 
\subsection*{Subsection $6.2'$. Splicing a knot complement with the complement of 
a trefoil}\label{subsec:splicing-with-T}

For a knot $K\subset Y$, let $Y(R,K)$  denote the three-manifold
obtained by splicing the complement of ${K}\subset {Y}$ 
with the complements of the 
right handed  trefoil. We  study 
the rank $r_r(K)$  of $\ov\HFT(Y(R,K))$  in this subsection.
With the notation of Subsection~\ref{subsec:simplifications}, $r_r(K)=i(\Dd'(R,K))$.
 Replacing the block forms of (21) 
in $\Dd'(R,K)$, we find 
\begin{displaymath}
\Dd'(R,K)=\colvec{
0&0&0&0&0&I&0&\ovl{B}_1B_0&0&0\\
\ovl{B}_\infty B_1 &0&0&0&0&0&0&0&0&I \\
0&\ovl{B}_\infty B_1 &0&0&\ovl{B}_\infty A_1 &0&0&I&0&0 \\
\ovl{D}_\infty B_1 &0&I&0&0&0&0&0&0&0\\
0&\ovl{D}_\infty B_1& 0& I& \ovl{D}_\infty A_1&0&0&0&0&0\\
I&0&0&0&\ovl{B}_0B_\infty &0&0&0&0& \ovl{B}_0X\\
0&0&0&0&0&\ovl{B}_0B_\infty &0&0&0& 0\\
0&0&0&0&0&0& I& \ovl{D}_1B_0&0&0\\
0&0&0&0& \ovl{D}_0B_\infty &0&0&0&I&\ovl{D}_0X\\
0&0&0&0&0&\ovl{D}_0B_\infty &0&0& 0& I\\
},
\end{displaymath}
where $A_\bullet=A_\bullet(K), B_\bullet=B_\bullet(K), 
C_\bullet=C_\bullet(K),D_\bullet=D_\bullet(K),
\ovl{A}_\bullet=\ovl{A}_\bullet(K),\ovl{B}_\bullet=\ovl{B}_\bullet(K),
\ovl{C}_\bullet=\ovl{C}_\bullet(K),\ovl{D}_\bullet=\ovl{D}_\bullet(K)$
 and $X=X(K)$ for $\bullet\in\{0,1,\infty\}$.
Doing a series of cancellations that correspond to the identity matrices which appear 
as the 
$$(1,6), (3,8), (4,3), (5,4), (6,1), (8,7), (9,9) \ \text{and}\ (10,10)$$
entries in the above block presentation
we obtain the equivalent matrix
\begin{equation}\label{eq:rank-d2}
\begin{split}
R_r(K):&=
\left(\begin{array}{ccc}
0&\ovl{B}_0X\ovl{B}_\infty\\
\ovl{X}B_\infty \ovl{B}_1&\ovl{X} B_\infty \ovl{A}_1
+\ovl{D}_0X\ovl{B}_\infty,\\
\end{array}\right)
\end{split}
\end{equation}
Where $\ovl{X}=\ovl{X}(R)=\ovl{B}_\infty(R) B_1(R)\ovl{B}_0(R)$.
\begin{cor} \label{cor:splicing-with-R}
For a knot $K\subset Y$ let $Y(R,K)$ denote the three-manifold obtained 
by splicing the complement of $K$ and the complement of the trefoil. Then
\begin{equation}\label{eq:rank-h1}
\begin{split}
\ov{\HFT}(Y(R,K)) &=\imath\big(R_r(K)\big).
\end{split}
\end{equation} 
\end{cor}
\begin{proof}
The claim follows immediately from the above discussion.
\end{proof}
For the  trefoils, our computations imply that 
\begin{displaymath}
\begin{split}
&\ovl{X}(R)B_\infty(R)=X(R)\ovl{B}_\infty(R)
=\ovl{X}(L)B_\infty(L)=X(L)\ovl{B}_\infty(L)=0\\ 
\Rightarrow\ \ \ &R_r(R)=R_r(L)=0 \\
\Rightarrow\ \ \ &\left|\ov{\HFT}(Y(R,R))\right|=7\ \ \ \text{and}\ \ \ 
\left|\ov{\HFT}(Y(R,L))\right|=9.
\end{split}
\end{displaymath}
The above computations agree with the computations of Hedden and Levine 
\cite{Matt-Adam}.
\begin{cor}\label{cor:splicing-with-trefoil}
For every knot $K$ in a homology sphere $Y$ we have
\begin{displaymath}
\begin{split}
\left|\ov{\HFT}(Y(R,K))\right|&\geq
 \big(a_0(K)+a_1(K)+2a_\infty(K)\big)-4\min\big\{a_0(K),a_1(K),a_\infty(K)\big\}\\
 &=4\max\big\{h_0(K),h_1(K),h_\infty(K)\big\}-\big(h_0(K)+h_1(K)+2h_\infty(K)\big).
\end{split}
\end{displaymath}
Moreover, if $Y(R,K)$ is a homology sphere $L$-space $K$ is trivial and $Y$ is a 
homology sphere $L$-space.
\end{cor}
\begin{proof}
Let $M=M(K)=X(K)\ovl{B}_\infty(K)$ and $\ovl{M}=\ovl{M}(K)
=\ovl{X}(K){B}_\infty(K)$ and note that 
\begin{displaymath}
\begin{split}
\rank(R_r(K))&=\rank\left(\begin{array}{cc}
0& \ovl{B}_0(K)M\\ \ovl{M}B_1(K)& \ovl{M}A_1(K)+\ovl{D}_0(K)M
\end{array}\right) \\
&\leq \rank\left(\begin{array}{cc}
 \ovl{M}B_1(K)& \ovl{M}A_1(K)
\end{array}\right) +\rank\left(\begin{array}{c}
\ovl{B}_0(K)M\\ \ovl{D}_0(K)M
\end{array}\right)\\
&=\rank(M)+\rank(\ovl{M}). 
\end{split}
\end{displaymath}
For every knot $K\subset Y$ as above note that the ranks of $M=M(K)$ 
and $\ovl{M}=\ovl{M}(K)$ are at most 
equal to the minimum of the sizes of the matrices $B_0(K), B_1(K)$ and $B_\infty(K)$,
which is 
$$\min\big\{a_0(K),a_1(K),a_\infty(K)\big\}.$$
Since $R_r(K)$ is of size $h_0(K)\times h_1(K)=(a_1(K)+a_\infty(K))\times (a_0(K)+a_\infty(K)$ 
this proves the first part of the corollary.\\

Let us assume that $\rank\left(\ov{\HFT}(Y(R,K))\right)=1$. From here we find
\begin{displaymath}
\begin{split}
&\big(a_0(K)+a_1(K)+2a_\infty(K)\big)-4\min\big\{a_0(K),a_1(K),a_\infty(K)\big\}\\
&\ \ \ \ \ \ \ \ \ \ \ \ \ \ \ \ 
=\big(a_0(K)+a_1(K)+2a_\infty(K)\big)-2\rank(M)-2\rank(\ovl{M})=1.
\end{split}
\end{displaymath}
Since $a_1(K)$ and $a_\infty(K)$ have the same parity while the parity of $a_0(K)$ is different 
from the parity of both $a_1(K)$ and $a_\infty(K)$, one can easily conclude that
 $a_0(K)-1=a_1(K)=a_\infty(K)$. Let $a$ denote the common value $a_1(K)=a_\infty(K)$. 
 Then the ranks of $M$ and $\ovl{M}$ are $a$
and $B_0(K), X(K)$ and $\ovl{X}(K)$ are all invertible. 
We may thus assume that $A_0(K)=D_0(K)=0$.
Since
\begin{displaymath}
\rank(\pphi_\infty(K)+\pphibar_\infty(K))=\rank \colvec{
B_1(K)A_0(K)& B_1(K)B_0(K)\\ I+D_1(K)A_0(K)& D_1(K)B_0(K)}=2a
\end{displaymath}
the three-manifold $Y$ is an $L$-space. Since splicing $K$ with the trefoil
is also a homology sphere $L$-space we conclude that $K$ is trivial, by 
Theorem 1 from  \cite{Matt-Adam}.
\end{proof}

\setcounter{thm}{0}
\def\thesection{\Alph{section}}
\setcounter{theorem}{0}
\setcounter{section}{1}  
\setcounter{equation}{23} 
\subsection*{Appendix$'$; Bordered Floer homology for knot complements}
The first draft of this paper appeared while the theory of bordered Floer homology 
was being developed. With bordered Floer homology conventions widely known to the 
Heegaard Floer community, the referee recommended the inclusion of an appendix 
which addresses the contribution of this paper within the realm of bordered Floer homology.\\

Let $K\subset Y$ denote a null-homologous knot inside the three-manifold $Y$, and 
let $H=(\Sig,\alphas,\ov\betas \cup\{\lambda,\mu\};z)$ denote a special Heegaard 
diagram for $K$, as constructed in Lemma~4.1. 
In particular, $H$ is a nice Heegaard diagram  for the bordered 
three-manifold $Y_K$ determined by $K\subset Y$ in the sense of \cite{LOT}. 
The Bordered Floer 
complex $\ov\CFDT(Y_K)$ may then be constructed from the chain complexes 
$M=M(K)$ and $L=L(K)$
(which are described in Proposition~5.1 as the mapping 
cones of $\pphibar_\infty(K)\co C_0(K)\ra C_1(K)$ and 
$\pphi_0(K)\co C_1(K)\ra C_\infty(K)$, respectively) and the chain maps 
$\Phi=\Phi(K)\co L\ra M$ and $\Psi_{i}=\Psi_i(K)\co M\ra L, i=1,2,3$.\\

More precisely and following the notation of Subsection~4.2 from 
\cite{LOT-E}, the idempotents
$\imath_0$ and $\imath_1$, and the chords $\rho_1,\rho_2,\rho_3,\rho_{12}=\rho_1\rho_2,
\rho_{23}=\rho_2\rho_3$ and $\rho_{123}=\rho_1\rho_2\rho_3$ 
form a $\Fbb$-basis for the differential graded 
algebra associated with the torus boundary;
\begin{displaymath}
\begin{split}&\\
&\mathcal{A}(T^2,0)=
\Big\langle\begin{diagram}[w=3em]
\imath_0 \ \bullet\ & \upperarrow{\rho_1}
\lift{-2}
{\strarrow{\rho_2}}
\lowerarrow{\rho_3} &\  \bullet\ \imath_1
\end{diagram}\Big\rangle/\left(\rho_2\rho_1=\rho_3\rho_2=0\right).\\
&
\end{split}
\end{displaymath}  
The module $\ov\CFDT(Y_K)$ is generated (over $\mathcal{A}(T^2,0)$) by 
the generators of $M$ and $L$. For  a generator $\x$ of $L$ we have
\begin{equation}\label{eq:d_L}
\begin{split}
& I(\x)=\imath_0\ \ \ \ \text{and}\ \ \ \
\partial (\x)=d_L(\x)+\rho_1\Psi_1(\x)+\rho_3\Psi_2(\x)+\rho_{123}\Psi_3(\x),
\end{split}
\end{equation}
while for a generator $\y$ of $M$ we have 
\begin{equation}\label{eq:d_M}
\begin{split}
& I(\y)=\imath_1\ \ \ \ \ \text{and}\ \ \ \ \ 
\partial (\y)=d_M(\y)+\rho_2\Phi(\y).
\end{split}
\end{equation}

The splicing formula of (17) is then just the gluing formula for bordered Floer homology, i.e. Theorem 1.3 from \cite{LOT}. A related discussion is 
carried over in Section 8 of \cite{LOT}.\\

\begin{defn}\label{def:admissible-data}
The chain complexes $(C_\bullet(K),d_\bullet), \bullet\in\{0,1,\infty\}$ and 
the chain maps $f_\bullet(K),\fbar_\bullet(K), \bullet\in\{0,\infty\}$   
are called {\emph{admissible data}} associated with the knot $K$ if they satisfy the following 
conditions:
\begin{itemize}
\item The homology of the complex $(C_\bullet(K),d_\bullet)$ is $\Hbb_\bullet(K)$.
\item The maps induced by $f_\bullet(K)$ and $\fbar_\bullet(K)$ in homology (under the 
identification of  the homology of  $(C_\bullet(K),d_\bullet)$ with $\Hbb_\bullet(K)$)
are $\pphi_\bullet(K)$ and $\pphibar_\bullet(K)$, respectively.
\item We have $f_0(K)\circ f_\infty(K)=\fbar_0(K)\circ \fbar_\infty(K)=0$.
\item The corresponding maps 
$$\theta(K)\co \Ker(\pphi_\infty(K))\ra \Coker(\pphi_0(K))\ \ \text{and}\ \ 
\thetabar(K)\co \Ker(\pphibar_\infty(K))\ra \Coker(\pphibar_0(K))$$
are isomorphisms, and are the inverses of the maps induces by $\pphi_1(K)$ and 
$\pphibar_1(K)$, respectively.
\end{itemize}
\end{defn}

 The proof of Theorem~1.1
implies that $(C_\bullet^i,d_\bullet^i)$ and 
the chain maps $f_\bullet^i,\fbar_\bullet^i$ for  $\bullet\in\{0,\infty\},i=1,2$
in (17) may be replaced by other admissible data corresponding 
to the knots $K_1$ and $K_2$. Correspondingly, the bordered Floer complex associated 
with any knot $K\subset Y$ may be constructed from admissible data associated with $K$.
More precisely, we have the following proposition.

\begin{prop}\label{prop:BFH}
Suppose that the chain complexes $(C_\bullet(K),d_\bullet), \bullet\in\{0,1,\infty\}$ and  
the chain maps $f_\bullet=f_\bullet(K),\fbar_\bullet=\fbar_\bullet(K), \bullet\in\{0,\infty\}$ 
are admissible data associated with the knot $K\subset Y$ and set
 \begin{displaymath}
 \begin{split}
 &M(K)=C_0(K)\oplus C_1(K),\ \ L(K)=C_1(K)\oplus C_\infty(K)\\
 \end{split}
 \end{displaymath}
 The bordered Floer complex $\ov\CFDT(Y_K)$ may then be constructed as the 
 left module over the differential graded algebra $\mathcal{A}(T^2,0)$ which is generated 
 by $\imath_0.L(K)$ and $\imath_1. M(K)$, and equipped with the differential 
 $\partial\co \ov\CFDT(Y_K)\ra \ov\CFDT(Y_K)$ defined by 
 \begin{equation}\label{eq:differential}
\partial \colvec{
\x\\ \y }=
\begin{cases}
\colvec{
d_0(\x)\\ \fbar_\infty(\x)+d_1(\y)}
+\colvec{
\rho_1f_\infty(\x)\\ \rho_3\fbar_0(\y)+
\rho_1\rho_2\rho_3\fbar_0(f_\infty(\x))}\ \ \ \ &\text{if }\colvec{\x\\ \y}\in M(K)\\
\colvec{
d_1(\x)\\ f_0(\x)+d_\infty(\y)}
+\rho_2.\colvec{
0\\ \x}
\ \ \ \ &\text{if }\colvec{\x\\ \y }\in L(K)\\
\end{cases} 
 \end{equation}
 \end{prop}

In particular, let the $\Fbb$-modules $\A_\bullet=\A_\bullet(K),\bullet\in\{0,1,\infty\}$ 
and the matrices 
$A_\bullet=A_\bullet(K),B_\bullet=B_\bullet(K),C_\bullet=C_\bullet(K)$ and 
$D_\bullet=D_\bullet(K)$ be defined as in 
Subsection~3.3. Set 
\begin{displaymath}
\begin{split}
&(C_0(K),d_0)=\left(\A_\infty\oplus \A_1,0\right),\ \ \ \ 
(C_\infty(K),d_\infty)=\left(\A_1\oplus \A_0,0\right)\\
&C_1(K)=\A_1\oplus \A_0\oplus \A_\infty\oplus \A_1\ \ \ \text{and}\ \ \ 
d_1=\left(\begin{array}{cccc}
0&0&0&0\\ 0&0&0&0\\ 0&0&0&0\\ I_{\A_1}&0&0&0\\ 
\end{array}\right).
\end{split}
\end{displaymath}
Correspondingly, define
\begin{displaymath}
\begin{split}
&f_\infty(K)=\left(\begin{array}{cc}
0&0\\ 0&0\\ I&0\\ 0&I
\end{array}\right),\ \ \ \ \ \ \ \ 
f_0(K)=\left(\begin{array}{cccc}
I&0&0&0\\ 0&I&0&0
\end{array}\right),\\ 
&\tau_1(K)=\left(\begin{array}{cccc}
0&0&0&0\\ 0&A_1&B_1&0\\ 0&C_1&D_1&0\\ 0&0&0&0
\end{array}\right),\ \ \ \ \ 
\ovl{\tau}_1(K)=\left(\begin{array}{cccc}
0&0&0&0\\ 0&\ovl{A}_1&\ovl{B}_1&0\\ 0&\ovl{C}_1&\ovl{D}_1&0\\ 0&0&0&0
\end{array}\right)
\end{split}
\end{displaymath}
and set $\fbar_\infty(K)=\tau_1(K) f_\infty(K)\tau_0^{-1}(K)$ and 
$\fbar_0(K)=\tau_\infty(K) f_0(K)\ovl{\tau}_1(K)$.
The data associated with $K$ 
consisting of $(C_\bullet(K),d_\bullet)$ and 
 $f_\bullet(K),\fbar_\bullet(K), \bullet\in\{0,\infty\}$ is then admissible. \\
 
 Corresponding to the above admissible data and associated with $K\subset Y$ 
  we may construct
 the bordered Floer complex for $K$ via
 \begin{displaymath}
 \begin{split}
 &M(K)=C_0(K)\oplus C_1(K)=
 \A_\infty\oplus \A_1\oplus \A_1\oplus \A_0\oplus \A_\infty\oplus \A_1\\
&L(K)=C_1(K)\oplus C_\infty(K)=
 \A_1\oplus \A_0\oplus \A_\infty\oplus \A_1\oplus \A_1\oplus \A_0\\ 
 &d_M=\left(\begin{array}{cccccc}
0&0&0&0&0&0\\
 0&0&0&0&0&0\\
 0&0&0&0&0&0\\
 B_1\ovl{A}_0& B_1\ovl{B}_0&0&0&0&0\\
D_1\ovl{A}_0& D_1\ovl{B}_0&I&0&0&0\\
 0&0&0&0&0&0 
 \end{array}\right),\ \ \ \ \ \ \ 
 d_L=\left(\begin{array}{cccccc}
 0&0&0&0&0&0\\
 0&0&0&0&0&0\\
 0&0&0&0&0&0\\
 I&0&0&0&0& 0\\
 0&0&I&0&0& 0\\
 0&0&0&I&0&0
 \end{array}\right)\\ 
 \end{split}
 \end{displaymath}
  \begin{displaymath}
 \begin{split}
 &\Phi(K)=\left(\begin{array}{cccccc}
 0&0&0&0&0&0\\
 0&0&0&0&0&0\\
 I&0&0&0&0&0\\
 0&I&0&0&0& 0\\
 0&0&I&0&0& 0\\
 0&0&0&I&0&0
 \end{array}\right),\ \ \ \ \ \ \ \ \ \ \ \Psi_1(K)=\left(\begin{array}{cccccc}
 0&0&0&0&0&0\\
 0&0&0&0&0&0\\
 I&0&0&0&0&0\\
 0&I&0&0&0& 0\\
 0&0&0&0&0& 0\\
 0&0&0&0&0&0
 \end{array}\right)\\
 &\Psi_2(K)=\left(\begin{array}{cccccc}
 0&0&0&0&0&0\\
 0&0&0&0&0&0\\
 0&0&0&0&0&0\\
 0&0&0&0&0& 0\\
 0&0&0&B_\infty \ovl{A}_1&B_\infty \ovl{B}_1& 0\\
 0&0&0&D_\infty \ovl{A}_1&D_\infty \ovl{B}_1&0
 \end{array}\right)\ \ \text{and}\ \ \Psi_3(K)=\Psi_2(K)\Phi(K)\Psi_1(K).
 \end{split}
 \end{displaymath}
 as the left module over the differential graded algebra $\mathcal{A}(T^2,0)$
 generated by $\imath_0.L$ and $\imath_1.M$
 and equipped with the differential 
 $\partial\co \ov\CFDT(Y_K)\ra \ov\CFDT(Y_K)$ defined by 
 the equations~(\ref{eq:d_L}) and (\ref{eq:d_M}).\\
 
 \begin{remark}
 Simultaneous computation of the matrices $\tau_\bullet(K)=\left(\begin{array}{cc}
A_\bullet&B_\bullet\\ C_\bullet& D_\bullet  \end{array}\right)$ is {\emph{a priori}} 
quite difficult,
as we observed in the case of trefoils in Section~6. This 
makes the above description of the bordered Floer homology hard to use even for knots 
$K\subset Y$ where we have complete understanding of the Heegaard Floer complex 
associated with $K$. However, it is possible to construct admissible data associated 
with $K\subset Y$  completely
in terms of the filtered chain complex 
$\CFT^\infty(Y,K;\Fbb)$, as will be discussed in the revision of 
\cite{Ef-essential}.  
 \end{remark}


\end{document}